\documentclass[12pt]{amsart}

\usepackage{amssymb,amscd,fullpage}

\newcommand{\p}{{\mathfrak p}}
\newcommand{\q}{{\mathfrak q}}
\newcommand{\oo}{{\mathcal O}}
\newcommand{\C}{{\mathbf C}}
\newcommand{\F}{{\mathbf F}}
\newcommand{\Q}{{\mathbf Q}}
\newcommand{\R}{{\mathbf R}}
\newcommand{\Z}{{\mathbf Z}}
\newcommand{\aaa}{{\mathfrak a}}

\newcommand{\vsp}{\vspace{8pt}}

\newcommand{\rarr}{{\rightarrow}}

\newcommand{\mymod}[2]{{ #1 \: (\bmod \: {#2})}}
\newcommand{\gln}{{GL_n}}

\newtheorem{thm}{Theorem}
\newtheorem{lem}{Lemma}
\newtheorem{prop}{Proposition}

\newtheorem*{defn}{Definition}

\theoremstyle{remark}
\newtheorem{remark}{Remark}

\newcommand\gal{{ \mbox{\rm Gal}}}
\newcommand{\ov}[1]{{\overline{{#1}}}}

\newcount\timehh\newcount\timemm
\timehh=\time 
\divide\timehh by 60 \timemm=\time
\begin{document}

\title[`Abstract' homomorphisms of $l$-adic groups]
	{`Abstract' homomorphisms of $l$-adic\\
		 groups and Abelian varieties}

\author{Siman Wong
  {
    \protect \protect\sc\today\ -- 
    \ifnum\timehh<10 0\fi\number\timehh\,:\,\ifnum\timemm<10 0\fi\number\timemm
    \protect \, \, \protect \bf DRAFT
  }
}

\address{Department of Mathematics \& Statistics, University of Massachusetts.
	Amherst, MA 01003-4515 USA}

\email{siman@math.umass.edu}

% \thanks{}


%\date{\today}


\begin{abstract}
Let $k$ be a totally real field, and let $A/k$ be an absolutely irreducible,
polarized Abelian variety of odd, prime dimension whose  endomorphism rings
is non-trivial and is defined over $k$.  Then the only
strictly compatible families of abstract, absolutely irreducible
representations
of $\gal(\ov{k}/k)$ coming from $A$ are tensor products of Tate twists of
symmetric powers of two-dimensional $\lambda$-adic
representations plus field automorphisms.
The main ingredients of the proofs are the work of Borel and Tits on the
`abstract' homomorphisms of almost simple algebraic groups, plus the work of
Shimura on the fields of moduli of Abelian varieties.
\end{abstract}

\maketitle

% \tableofcontents

\section{Introduction}

Let $A$ be an absolutely irreducible Abelian threefold over a number field $k$.
Suppose that the
endomorphisms of $A$ are all defined over $k$.  Let $l$ be a prime, and write
$
T_l(A)
$
for the $l$-adic Tate module of $A$.  Let $L/\Q_l$ be a finite extension,
set
$
G_k := \gal(\ov{k}/k)
$,
and let
$
\rho: G_k \rarr GL_3(L)
$
be an \textit{abstract} (i.e.~not necessarily continuous) group homomorphism
which is a quotient of the $l$-adic
representation
$
G_k \rarr \text{Aut}(T_l(A))
$
by an abstract homomorphism
$
\text{Aut}(T_l(A)) \rarr GL_3(L)
$.
A special case of our main results (Theorems \ref{thm:main} and
\ref{thm:indep}) says that
if $l$ is sufficiently large and if $k$ is \textit{totally real}, then
$\rho$ is the Tate
twist of the symmetric-square of an $l$-adic representation 
$
G_k \rarr GL_2(L)
$
plus a field automorphism $L\rarr L$.   To put this result in perspective, we
first recall
how $l$-adic representations arise in the Langlands program.

\vsp

Denote by $\mathbb A_k$ the adeles of a number field $k$.  The Langlands
program predicts that algebraic cuspidal automorphic forms on
$
GL_n(\mathbb A_k)
$
correspond to motives, or, more concretely, compatible families of
$n$-dimension $l$-adic
representations coming from \'etale cohomology of algebraic varieties
defined over $k$.  For
\textit{holomorphic} eigenforms on
$
GL_2(\mathbb A_\Q)
$
this conjecture follows from the work of Eichler-Shimura, Deligne, and
Deligne-Serre, 
by relating the eigenforms to the cohomology of certain Kuga fiber varieties.
This has been extended to holomorphic $GL_2$-eigenforms over totally real
fields,
by Blasius-Rogawski  \cite{blasius} and Taylor \cite{taylor}.  When the
automorphic form corresponds to a selfdual cuspidal representation,
Clozel establishes the existence of these $l$-adic representations
under mild conditions \cite[p.~150]{clozel}.  Beyond that little is known
for non-selfdual representations for $n\ge 3$.

\vsp

Clozel raises the question
of whether we can collect numerical evidence to support the Langlands
correspondence for non-selfdual representations.
Since then several authors have carried out numerical studies in the
smallest open case $n=3$
	(\cite{ash2}, \cite{ash3}, \cite{top1}, \cite{top2}).
Such investigations have two parts.  On the automorphic side,
there is a well-known connection between cuspidal automorphic forms on
$
GL_3(\mathbb A_\Q)
$
and cuspidal cohomology of arithmetic subgroup of 
$
SL_3(\Z)
$
(cf.~\cite{schwermer} for an exposition).  Thus automorphic forms on $GL(3)$
can be systematically enumerated, at least for small levels, via Ash's modular
symbol algorithm on $GL(3)$ \cite{ash1}.  On the Galois side, however, there
is no
known algorithm that generates all irreducible $3$-dimensional $l$-adic
representations.  Examples in the literature come from either ingenious but
ad hoc geometric constructions (\cite{top1}, \cite{top2}), or, in the case of
mod-$l$ representations, number fields with special Galois groups
	(\cite{ash2}, \cite{ash3}).

\vsp

A familiar source of compatible families of $l$-adic representations comes from
the Tate modules of Abelian varieties.  It is natural
to try to extract from them
non-selfdual representations.
To be specific, denote by $F_\lambda$ the completion of a number field $F$ at
a finite place
$\lambda$ of residual characteristic $l$.  Suppose an abstract
homomorphism\footnote{we
	will explain shortly will we want to work with not-necessarily
	continuous maps.}
$
\rho_\lambda: G_k \rarr GL_N(F_\lambda)
$
\textit{comes from} $A/k$, i.e.~there exists an abstract homomorphism
$
\iota_\lambda: \text{Aut}(T_l(A)) \rarr GL_N(\ov{F_\lambda})
$
such that, along with the canonical inclusion
$$
j_\lambda: GL_N(F_\lambda) \hookrightarrow GL_N(\ov{F_\lambda}),
$$
the following diagram commutes ($V_l(A) := T_l(A)\otimes\Q_l$):
\if 3\
{
\begin{equation}
\divide\dgARROWLENGTH by2
\begin{diagram}
\node{G_k}
	\arrow[2]{e,t}{\rho_{A, l}}	\arrow{s,l}{\rho_\lambda}
\node[2]{\text{Aut}(T_l(A))}
\node{\subset}
\node{\text{Aut}(V_l(A))}
	\arrow{s,r}{\iota_\lambda}
\\
\node{GL_N(F_\lambda)}
	\arrow[4]{e,t}{j_\lambda}
\node[4]{GL_N(\ov{F_\lambda}).}
\end{diagram}
	\label{diagram}
\end{equation}
}
\fi

\noindent
\begin{minipage}[b]{.2in}
	\begin{equation}
	\label{diagram}
	\end{equation}

\mbox{}

\mbox{}
\end{minipage}
\begin{minipage}[t]{6in}
\begin{picture}(400,100)

\put(100, 10){$GL_N(F_\lambda)$}
\put(285, 10){$GL_N(\ov{F_\lambda})$}
\put(210, 19){$\scriptstyle j_\lambda$}

\put(155, 13){\vector(1,0){120}}
\put(105, 51){$\scriptstyle \rho_\lambda$}

\put(113, 77){$G_k$}
\put(200, 77){$\text{Aut}(T_l(A)) \:\: \subset \:\: \text{Aut}(V_l(A))$}
\put(305, 67){\vector(0,-1){40}}

\put(133, 80){\vector(1,0){60}}
\put(150, 86){$\scriptstyle \rho_{A, l}$}

\put(120, 67){\vector(0,-1){40}}
\put(313, 50){$\scriptstyle \iota_\lambda$}
\end{picture}
\end{minipage}

\if 3\
{
\begin{equation}
\begin{myCD}
G_k @>{\rho_{A, l}}>>	\text{Aut}(T_l(A))\subset \text{Aut}(V_l(A))
\\
@V{\rho_\lambda}VV	@VV{\iota_\lambda}V
\\
GL_N(F_\lambda) @>{j_\lambda}>> GL_N(\ov{F_\lambda})
\end{myCD}
%	\label{diagram}
\end{equation}
}
\fi

Suppose $A/k$ is an Abelian threefold such that
$
\text{End}(A/k)\otimes\Q
$
contains a quadratic field.  Then
$V_l(A)$ splits into a direct sum of two conjugate, odd-dimensional Galois
representations.  Furthermore,  if
$
\text{End}(A/\ov{k}) \otimes\Q
$
is exactly a quadratic
field, by the Mumford-Tate conjecture
	\cite{serre:76}
 we expect  the Lie algebra of the
image of these conjugate representations to contain
$
{\mathfrak{sl}}_3(\Z_l)
$.
 The Jacobian of a generic Picard quartic
$
y^3 = x^4 + ax^3 + bx^2 + cx + d
$
over $k$ would satisfy the required conditions, thereby yielding interesting
irreducible,
non-selfdual $3$-dimensional $l$-adic representations, \textit{provided that}
$k$ contains a primitive third-root of unity.  One important consequence of
the result stated in the first paragraph
is that such examples cannot exist if $k$ is not totally real, even if we
do not require that the representation be continuous.

\vsp

The proofs of our main results have two parts.  As a preliminary step, we
utilize Shimura's
work on the fields of moduli of Abelian varieties with PEL structures to 
limit the possibilities for the endomorphism algebras
of odd-dimensional Abelian varieties over a \textit{totally real} field.

\vsp

\begin{thm}
	\label{thm:real}
Let $A$ be an absolutely simple, polarized Abelian variety of odd
dimension $d$ defined over a totally real number field $k$.  Suppose that the
endomorphisms of $A$ are all defined over $k$.   Then
$
E := \text{End}(A/k)\otimes \Q
$
is a totally real number field of degree dividing $d$.
\end{thm}

\vsp

With $A/k$ as in Theorem \ref{thm:real}, set $[E:\Q]=g$ and 
$
\delta = d/g
$.
Since the action of $E$ on $V_l(A)$ commutes that of $G_k$, we have a direct
sum decomposition
\begin{equation}
\rho_{A, l} \simeq \tau_{l, 1}\oplus \cdots \oplus\tau_{l, g},
	\label{direct}
\end{equation}
where each $\tau_{l, m}$ is a representation
$
G_k\rarr GL_{2\delta}(E_{l, m})
$
for some completion $E_{l, m}$ of $E$ at a place above $l$.  Denote by
$$
\pi_{l, m}: \text{Aut}(T_l(A)) \rarr GL_{2\delta}(E_{l, m})
$$
the corresponding quotient.   If
$
\rho_\lambda
$
is irreducible, it must factor through one of the $\tau_{l, m}$.  If
$
\iota_\lambda
$
is further assumed to be $l$-adic analytic, standard Lie theory then implies
that
$
\iota_\lambda
$
is the composition of some
$
\pi_{l, m}
$
with a rational\footnote{i.e.~a
	morphism of algebraic		\label{foot}
	groups.}
representation
$GL_{2\delta}(E_{l, m})\rarr GL_N(\ov{F_\lambda})$.  However, if we start out
with a $\lambda$-adic representations that arise in arithmetic geometry, we
only know that $\rho_\lambda$ is continuous.  So the additional hypothesis
that $\rho_\lambda$ comes from an Abelian variety does not automatically 
imply that $\iota_\lambda$ is continuous, let alone being $\lambda$-adic
analytic.  To address this issue, we invoke a deep result of 
Borel and Tits \cite{borel} which asserts that, over an infinite field $K$,
the `abstract' homomorphisms of an almost simple matrix group $G$ over $K$
into $PGL_n(\ov{K})$ are tensor products of the rational projective
representations of $G$
plus field automorphisms; additional work is required to handle ordinary
representations.
 Before we state our main results we introduce
some notation.

\vsp

Let $V$ be a two-dimension vector space over a field $K$.  Let
$
\varphi: K\rarr \ov{K}
$
be a field homomorphism.  For any integer $n\ge 1$, denote by
$
\sigma^{(n)}_\varphi \circ \varphi
$
the $n$-th symmetric power of the `$\varphi$-twisted' standard representation
$
GL_2(K)\rarr GL_2(\varphi(K))\subset GL_2(\ov{K})
$.
This is a $(n+1)$-dimensional representation.

\vsp

\begin{thm}
	\label{thm:main}
With the notation as in Theorem \ref{thm:real}, suppose further that
there exists an abstract, absolutely irreducible representation
$
\rho_\lambda: G_k \rarr GL_N(F_\lambda)
$
coming from $A$.  Suppose $l = \text{char}(\lambda)$ is
sufficiently large.

{\rm(a)} If $N < 2d$ and $\rho_\lambda$ is non-Abelian then $[E:\Q]>1$.

{\rm(b)} Suppose $d = [E:\Q]$.  Then there exist
\begin{itemize}
\item
a choice $\tau_{l, m}$ among the factors in the decomposition
	{\rm(\ref{direct})},
\item
positive integers $r(\lambda)$ and
$
n_1(\lambda), \ldots, n_{r(\lambda)}(\lambda)
$
with $\sum_{i=1}^{r(\lambda)} (n_i(\lambda)+1) = N$,
\item
a field homomorphism
$
\varphi_i(\lambda): E_{l, m}\rarr \ov{E_{l, m}}
$
for each $1\le i\le r(\lambda)$, and
\item
a homomorphism of multiplicative groups
$
f_\lambda: E_{l, m}^\times\rarr \ov{E_{l, m}}^\times
$,
\end{itemize}
such that, in the notation of diagram {\rm(\ref{diagram})}, 
$$
\iota_\lambda
=
\Bigl(
\bigl[
\otimes_{i=1}^{r(\lambda)}
	\:\:
  \sigma^{(n_i(\lambda))}_{\varphi_i(\lambda)}\circ\varphi_i(\lambda)
\bigr]
\otimes (f_\lambda\circ\det)
\Bigr)
\circ	\tau_{l, m}.
$$
\end{thm}

\vsp

\begin{remark}
(a) Set $d=N=3$ and we recover the result stated in the first
paragraph.

(b)
Our technique applies \textit{mutatis mutandis} when $[E:\Q] < d$,
once we replace the symmetric power representations of $GL_2$ by the Weyl
modules associated to $GL_n$
	\cite[$\S 15.5$]{fulton-harris}.
The same goes for Theorem \ref{thm:indep} below.
In addition, Theorems \ref{thm:main} and \ref{thm:indep} remain true if we
replace $\lambda$-adic
representations by mod-$\lambda$ representations.

\end{remark}

\vsp

In practice, continuous $\lambda$-adic representations often arise in a
compatible family.  Our next result says that if
$
\{ \rho_\lambda \}_\lambda
$
is a compatible family coming from an Abelian variety $A/k$ as in Theorem
\ref{thm:real}, then the parameters
$
r(\lambda), n_i(\lambda), f_\lambda
$
and $\varphi_i(\lambda)$ in Theorem \ref{thm:main} are essentially independent
of $\lambda$.  To set up the notation we first
we review the definition of compatible families of representations.

\vsp

\begin{defn}
A {\rm strictly compatible family of $\lambda$-adic rational representations of $G_k$}
consists of
\begin{itemize}
\item
a number field $F/\Q$,
\item
a collection
$
\{ \rho_\lambda: G_k \rarr \gln(F_\lambda) \}_\lambda
$
of continuous Galois representations for each finite place $\lambda$ of $F$,
\item
a finite set $S$ of finite places of $\oo_k$,
\end{itemize}
such that, with
$
\Phi_\lambda: F\rarr F_\lambda
$
denoting the canonical embedding,
\begin{itemize}
\item
for any prime $\p\not\in S\cup \{ \lambda \}$, the coefficients of the
	characteristic polynomial
$
\text{char}(\rho_\lambda, \p)
$
of
$
\rho_\lambda(Frob_\p)
$
lies in $\Phi_\lambda(F)$, and
\item
for any two distinct finite places $\lambda, \q$ and for any place
$
\p\not\in S\cup\{\lambda, \q\}
$,
$$
\Phi_\lambda^{-1}(\text{char}(\rho_\lambda, \p))
=
\Phi_\q^{-1}(\text{char}(\rho_\q, \p)).
$$
\end{itemize}
\end{defn}

\vsp

\if 3\
{
(a) 
Note that if $\{ \rho_\lambda \}_\lambda$ is a compatible family, then
$\rho_\lambda$ is
irreducible for \textit{every} $l$ if and only if $\rho_\lambda$ is
irreducible for
\textit{one} $\lambda$.

(b)
}
\fi

\begin{remark}
The phrase `rational representations'  refers to the condition
that the image of each $\rho_\lambda$ lies in $\Phi_\lambda(F)$ for a fixed
number
field $F$, rather than the whole $\ov{F_\lambda}$.  In particular, this is
different from the terminology in the theory of algebraic groups
(cf.~footnote \ref{foot}).
\end{remark}

\vsp

\begin{thm}
	\label{thm:indep}
Let $A/k$ be as in Theorem \ref{thm:real} with $[E:\Q]=d$.  Suppose
$
\{ \rho_\lambda \}_\lambda
$
is a strictly compatible family of absolutely irreducible rational representations coming
from $A/k$.  Following the notation in Theorem \ref{thm:main}, for all but
finitely many pairs
$
(\lambda, \q)
$
of finite places of $k$,
\newcounter{property}
\begin{list}
	{{\rm(\roman{property})}}{\usecounter{property}
				\setlength{\labelwidth}{30pt}}
\item
$r(\lambda) = r(\q)$,
\item
$
n_i(\lambda) = n_i(\q)
$
for all $i$, up to reordering,
\item
${f_\lambda}\bigr|_\Z = {f_\q}\bigr|_\Z$,  and
\item
there exists a field isomorphism
$
\xi_{\lambda, \q} : \ov{\Phi_\lambda(E)} \rarr \ov{\Phi_\q(E)}
$
such that
$$
\xi_{\lambda, \q}
	\bigl( \varphi_i(\lambda)\bigr|_{\Phi_\lambda(E)} \bigr)
=
	\bigl( \varphi_i(\q)\bigr|_{\Phi_\q(E)} \bigr)
\hspace{20pt}
\text{for every $1\le i\le r(\lambda)$.}
$$
\end{list}
\end{thm}

\vsp

\begin{remark}
We clarify the meaning of conditions (iii) and (iv).  In the notation of
Theorem \ref{thm:main} (and dropping the $\lambda$ from $r(\lambda)$ and
$n_i(\lambda)$),
$$
\iota_\lambda
=
\Bigl(
\bigl[
\otimes_{i=1}^r
  \sigma_{\varphi_i}^{(n_i)}
	\circ
  \varphi_i(\lambda)
\bigr]
\otimes (f_\lambda\circ\det)
\Bigr)
\circ
\tau_{l, m}.
$$
Since
$
\Phi_\lambda(E)\subset E_\lambda
$
is dense in $E_\lambda$, for the purpose of computing
$
\iota_\lambda
$,
condition (iv) completely determines
$
\varphi_\lambda
$.
Similarly, for $A/k$ as in Theorem \ref{thm:real},
$
\det \tau_{l, m}
$
takes values in $\Z_l-\{0\}$ and $\Z$ is dense in $\Z_l$, so condition (iii)
completely determines $f_\lambda$.
\end{remark}

\vsp

Since $G_k$ is compact, if
$
\rho_\lambda: G_k \rarr GL_N(F_\lambda)
$
is continuous then its image is conjugate to a subgroup of
$
GL_N(\oo_\lambda)
$,
where $\oo_\lambda$ denotes the ring of integers of $F_\lambda$
	(cf.~\cite[Remark 1 on p.~1]{serre:abelian}).
To say that $\rho_\lambda$ comes from $A/k$ is to say that the
	$\oo_\lambda[G_k]$-module
for $\rho_\lambda$,
	\textit{after tensoring with $\Q$},
is a quotient of $V_l(A)$.  This suggests that we seek an \textit{integral}
form of our results,
i.e.~\textit{without} taking the tensor product.  This calls for an
integral form of the work of Borel-Tits, i.e.~a
characterization of abstract homomorphisms of linear groups such as
$
GL_n(\oo_\lambda)
$
into
$
GL_m(\oo'_\lambda)
$.
This is a very interesting question by itself, but little is known in this
direction.  For the record we have the following partial results.

\vsp

\begin{thm}
	\label{thm:integral}
For any integers $n, N>1$ there exists a constant $c(n, N)>0$ such that, for
any prime
$
l > c(n, N)
$
and any two rings of integers $\oo$ and $\oo'$ of finite extensions of
$\Q_l$,

\mbox{\rm (a)}
every abstract homomorphism
$
\varphi: SL_n(\oo')\rarr GL_N(\oo)
$
is continuous and $\ker \varphi$ is contained in the subgroup of scalar
matrices; and

\mbox{\rm (b)}
if $N < n$, then there is no non-trivial abstract homomorphism
$
\varphi: SL_n(\oo') \rarr GL_N(\oo)
$.

\noindent
Part {\rm(a)}
 remains true if we replace $SL_n$ by $GL_n, Sp_{2n}, PSL_n, PGL_n$, or
$
PSp_{2n}
$.
The same goes for Part {\rm(b)}, if for $GL_n$ and $PGL_n$ we require that
$\varphi$ be
non-Abelian;
and for $Sp_{2n}$ and $PSp_{2n}$, that $N<2n$.

\end{thm}

\section{Moduli of Abelian varieties}
	\label{sec:moduli}

In this section we recall results of Shimura on Abelian varieties with
prescribed endomorphism structures.
Let $A$ be an Abelian variety over $\C$. Let $E$ be a semisimple algebra over
a field $k$.  Let
$
\theta: E\rightarrow End_\C(A)\otimes_\Z \Q
$
be an isomorphism,
and let $\mathcal C$ be a polarization of $A$. Two such triples
$
(A, \mathcal C, \theta)
$
and
$
(A', \mathcal C', \theta')
$
are said to be isomorphic to each other if there exists an isomorphism
$
A\stackrel{\lambda}{\rarr}A'
$
such that
$
\lambda(X')\in\mathcal C
$
for every divisor $X'\in\mathcal C'$, and that
$
\lambda\theta(\alpha) = \theta'(\alpha)\lambda
$
for every $\alpha\in E$.

\begin{prop}[Shimura {\cite[Prop.~5]{shimura_autofcn}}]
	\label{prop:moduli}
With
$
\mathcal P = (A, \mathcal C, \theta)
$
as before, suppose further that $A$ is defined over a number field. Then
there exists a unique subfield $k/\Q$, called the
{\rm field of moduli} of $(A, \mathcal C, \theta)$, which is
 characterized by the following
property:
For every automorphism $\sigma$ of $\C$, the triple $\mathcal P^\sigma$ is
isomorphic to $\mathcal P$ if and only if
$
\bigl.\sigma\bigr|_k = id$.

	\qed
\end{prop}

Let $F$ be a totally real
number field contained in the center of $E$. Let $\rho$ be an involution of
$E/k$ such that $a^\rho = a$ for every $a\in F$, and that
$
tr_{E/\Q}(a a^\rho)>0
$
for every non-zero $a\in E$. Let $n=\dim A$, and let $\Phi$ be a
representation of $L$ by
the $n\times n$ complex matrices containing no zero representations.
Following Shimura \cite[$\S 1.4$]{shimura_fielddefn}, the triple
$
(A, \mathcal C, \theta)
$
is called a \textit{polarized Abelian variety of type} $\{E, \Phi, \rho\}$
if
\begin{list}
	{(\roman{property})}{\usecounter{property}
				\setlength{\labelwidth}{30pt}}
\item
$\Phi$ takes the identity element of $E$ to the identity matrix.
\item
the representation of $\theta(x)$ for $x\in E$ with respect to any analytic
coordinates given by 
\\
$A\simeq \C^n/\text{(lattice)}$ is equivalent to $\Phi$.
\item
the involution of $End_\C(A)\otimes\Q$ determined by $\mathcal C$ coincides
on $\theta(E)$ with the involution 
\\
$\theta(a)\rarr\theta(a^\rho)$.
\end{list}

\begin{prop}[Shimura {\cite[Prop.~1.6]{shimura_fielddefn}}]
	\label{prop:trace}
Let $\mathcal (A, \mathcal C, \theta)$ be a polarized Abelian variety of
type $\{E, \Phi, \rho\}$. Then the field of moduli of the triple
$
(A, \mathcal C, \theta)
$
contains the number field generated over $\Q$ by the trace of the elements
$
\Phi(\alpha)
$
for all $\alpha$ in the center of $E$.

	\qed
\end{prop}

\vsp

Next, we recall Albert's work on the classification of endomorphism algebras
of simple complex Abelian varieties.  Let $A$ be a simple, polarized Abelian
variety over the complex numbers. Then
$
E := \text{End}(A)\otimes \Q
$
is a finite dimensional division algebra over $\Q$ equipped with the Rosati
involution with respect to the polarization.  Denote by $K$ the center of $E$,
and by
$K_0$ the fixed field of the Rosati involution restricted to $K$.  Set
$$
[E:K] = \delta^2,  [K:\Q]=e,  [K_0=\Q]=e_0.
$$
It is known that
$
\delta^2 e
$
divides $2\dim A$, and if $e=2\dim A$, then $K$ is a CM field
\cite[$\S 2.2$]{shimura_book}.

$E$ is said to be a division algebra of the first kind (resp.~second kind) if
the Rosati involution acts trivially on $K$ (resp.~non-trivally).
By the work of Albert, $E$ falls into one of four types, exactly three of
which correspond to division algebras of the first kind.  Moreover, we have
the following restrictions on $e, e_0$ and $\delta$
	(\cite[Prop.~5.5.7]{lange}):
\begin{table}[h]
\renewcommand{\arraystretch}{1.2}
\begin{tabular}{c|c|c|c}
End$(A)\otimes \Q$
&
$\delta$
&
$e_0$
&
restriction
\\	\hline
\text{totally real number field}
&
$1$ & $e$ & $e|\dim A$
\\	\hline
\text{totally indefinite quaternion algebra}
&
$2$ & $e$ & $2e|\dim A$
\\	\hline
\text{totally definite quaternion algebra}
&
$2$ & $e$ & $2e|\dim A$
\\	\hline
\text{division algebra of the second kind}
&
$\delta$ & $e/2$ & $e_0 \delta^2 |\dim A$
\end{tabular}
	\vspace{10pt}
	\caption{}
\end{table}

\renewcommand{\arraystretch}{1}

\section{Endomorphism algebras over totally real fields}

Let
$
(A, \mathcal C, \theta)
$
be an absolutely simple, polarized Abelian variety of type
$
\{ E, \Phi, \rho \}
$
over $k$ and of \textit{odd} dimension $d$.  Cf.~Table 1 and we see that
$
E = End_{\ov{k}}(A)\otimes\Q
$
is either a totally real number field of degree dividing $d$ or a division
algebra of the second kind.  To prove Theorem \ref{thm:real} we need to
eliminate this second possibility when $k$ is totally real.

Suppose $E$ is a division algebra of the second kind.  Following the notation
in Table I, we have
	\cite[$\S 2$]{shimura_analyticfamily}
\begin{equation}
E\otimes_\Q \R
\simeq
\underbrace{M_\delta(\C) \oplus \cdots \oplus M_\delta(\C)}_{e_0} ,
	\label{tensor}
\end{equation}
where $M_\delta(\C)$ denotes the algebra of
$
\delta\times\delta
$
complex matrices.  For
$
i=1, \ldots, e_0
$,
denote by
$
\chi_i: E\otimes_\Q \R\rarr M_\delta(\C)
$
the projection onto the $i$-th factor.  According to
	\cite[$\S 2$]{shimura_analyticfamily},
the inequivalent absolutely irreducible representations of $E$ are given by 
$$
\chi_1, \ov{\chi}_1, \cdots, \chi_g, \ov{\chi}_g,
$$
where
$
\ov{\chi}_i
$
denotes the complex conjugate of $\chi_i$.  Moreover, for each $i$ there
exist integers
$
r_i, s_i\ge 0
$
satisfying
\begin{eqnarray}
r_i + s_i
&=&
2\delta (\dim A)/[E:\Q]
	\nonumber
\\
&=&
(\dim A)/(\delta e_0)
	\hspace{40pt}
	\text{since $[E:\Q]=[E:K] [K:K_0] [K_0:\Q]$}
\label{ri}
\\
&&
	\hspace{178pt}
	= \delta^2 \cdot 2 \cdot e_0,
\nonumber
\end{eqnarray}
such that with
\begin{equation}
\Phi_i := (\chi_i \otimes I_{r_i}) \oplus (\ov{\chi}_i \otimes I_{s_i})
	\label{phii}
\end{equation}
($I_j$ denotes the $j\times j$ identity matrix), the type
$
\{ E, \Phi, \rho \}
$
of $A/k$ satisfies
	\cite[$\S 2$]{shimura_analyticfamily},
\begin{equation}
\Phi \simeq \Phi_1 \oplus \cdots \oplus \Phi_g.
	\label{phi}
\end{equation}
Apply this to the case where $\dim A$ is odd, it then follows from (\ref{ri})
that
$
r_i + s_i
$
is odd.  In particular,
$
r_i\neq  s_i
$
for every $i$.

In light of the isomorphism (\ref{tensor}), we can find $\epsilon\in E$ such
that
\begin{list}
	{(\roman{property})}{\usecounter{property}
				\setlength{\labelwidth}{15pt}}
\item
$\text{trace}(\chi_1(\epsilon))$ is complex, and
\item
the imaginary part of $\chi_1(\epsilon)$ is large, but
$
|\text{trace}(\chi_i(\epsilon))|
$
is small for every $i>1$.
\end{list}
Then
\begin{eqnarray}
\text{trace}(\Phi(\epsilon))
&=&
\sum_i \text{trace}(\Phi_i(\epsilon))
	\hspace{125pt}
	\text{by (\ref{phi})}
	\nonumber
\\
&=&
\sum_i
  \bigl[
    r_i \cdot \text{trace}(\chi_i(\epsilon)) + s_i \cdot \ov{\text{trace}(\chi_i(\epsilon))}
  \bigr]
	\hspace{10pt}
	\text{by (\ref{phii}).}
	\label{what}
\end{eqnarray}
Since
$
r_1\neq s_1
$,
condition (i) implies that the term $i=1$ in (\ref{what}) is not real.
Condition (ii) then implies that
$
\text{trace}(\Phi(\epsilon))
$
is not real.  By Proposition \ref{prop:trace}, the field of moduli of
$
(A, \mathcal C, \theta)
$
cannot be totally real.  On the other hand, by Proposition \ref{prop:moduli}
this field of moduli is contained in $k$, which is totally real by
hypothesis. This is a contradiction, and Theorem \ref{thm:real} follows.

\vsp

\section{`Abstract' homomorphisms of $l$-adic groups}
	\label{sec:borel-tits}

Let $G$ be a linear algebraic group (i.e.~a closed subgroup of $GL_n$)
defined over a field $K$.  Fix an algebraic closure $\ov{K}$ of $K$.
The group $G$ is said to be almost simple (resp.~absolutely almost simple) if
it has no proper, normal subgroup over $K$ (resp.~over $\ov{K}$) of dimension
$>0$.  It is said to be isotropic over $K$ if it
contains a non-zero, split $K$-torus.  Both $SL_n(K)$ and $Sp_{2n}(K)$ are
isotropic and are absolutely almost simple over $\ov{K}$.

Let $\varphi: K\rarr K'$ be a non-trivial field homomorphism.  Denote by
$
G_\varphi
$
the base change of $G$ to $K'$ via $\varphi$, i.e.~the pullback of
$
G\rarr \text{Spec } K
$
via $\varphi$:
$$
\begin{CD}
G_\varphi @>>> G
\\
@VVV @VVV
\\
\text{Spec } K' @>>> \text{Spec } K.
\end{CD}
$$
If 
$
\rho: G\rarr PGL_n(K')
$
is a rational representation, i.e.~a morphism of algebraic groups, then
we get a natural map
$
\rho_\varphi:  G_\varphi\rarr PGL_n(K')
$.

The following special case of \cite[Thm.~10.3]{borel} contains all the
ingredients from the work of Borel-Tits for our subsequent applications.

\begin{thm}[Borel-Tits]
	\label{thm:borel}
Let $K$ be a field of characteristics zero, and let $K'$ be an
algebraically closed field.  Let $G$ be a connected linear algebraic group
over $K$ which is absolutely almost simple and is isotropic over $K$.  Let
$
\rho: G(K) \rarr PGL_n(K')
$
be an irreducible projective representation {\rm(}$n\ge 2${\rm)}.  Then there
exists finitely many distinct, non-trivial field homomorphisms
$
\varphi_i: K\rarr K'
$,
and for each $i$, a non-trivial, irreducible rational projective
representation 
$
\pi(i)
$
of $G$, such that
$$
\rho \simeq   \otimes_i		(\pi(i)_{\varphi_i}\circ \varphi_i).
$$
In particular, char$(K')=0$.  The pairs consisting of $\varphi_i$ and the
equivalence class of $\pi_i$ are unique up to reordering.

	\qed
\end{thm}

The following fact is elementary.

\begin{lem}
	\label{lem:central}
Let
$
1\rarr C\rarr G\rarr \tilde{G} \rarr 1
$
be a central extension of groups.  Given any group $H$, any two lifts of a
homomorphism
$
H\rarr \tilde{G}
$
differs by the multiple of a homomorphism $H\rarr C$.

	\qed
\end{lem}

For the rest of this section, $K$ denotes a field of characteristic zero.
Let $i\not=j$ be integer between $1$ and $n$.  For any elements
$
\alpha\in K^\times
$
and
$
\beta\in K
$,
define two matrices
\begin{eqnarray*}
A_i(\alpha)
&=&
\text{the matrix obtained from the $2\times 2$ identity matrix by replacing the}
\\
&&
\text{$ii$-th entry by $\alpha$,}
\\
B_{ij}(\beta)
&=&
\text{the matrix obtained from the $2\times 2$ identity matrix by adding
$\beta$ to the}
\\
&&
\text{$ij$-th entry}.
\end{eqnarray*}
For $i\not=j$ we have the identity
\begin{equation}
A_i(\alpha) B_{ij}(\beta) A_i(\alpha)^{-1}
=
B_{ij}(\alpha^{\text{sgn}(j-i)}\beta).
	\label{matrix}
\end{equation}

\vsp

\begin{prop}
	\label{prop:gln}
Every abstract, absolutely irreducible representation
$
GL_2(K)\rarr GL_d(\ov{K})
$
is equivalent to 
$$
\bigl[
  \otimes_{s=1}^r ( \sigma_{\varphi_s}^{(n_s)} \circ {\varphi_s} )
\bigr]
  \otimes (f\circ\det),
$$
where $V$ is the underlying representation space associated to natural
representation of $GL_2(K)$;
 $\det: GL_2(K)\rarr K^\times$ is the determinant map; and
$
f: K^\times\rarr \ov{K}^\times
$
is a homomorphism of multiplicative groups.
\end{prop}

\begin{proof}
Since $\rho$ is absolutely irreducible, so are
$
\rho|_{SL_2(K)}
$
and
$$
\tilde{\rho}_S :=
\text{the projective representation $SL_2(K)\rarr PGL_d(\ov{K})$ associated to
$
\rho|_{SL_2(K)}
$}.
$$
Define
$$
\widetilde{\sigma_{\varphi_s}^{(n_s)}} :=
	\text{the projective representation associated to
	$\sigma_{\varphi_s}^{(n_s)}$}.
$$
The absolutely irreducible rational projective representations of $SL_2(K)$
are precisely the restrictions to $SL_2(K)$ of the
$
\widetilde{\sigma_{\varphi_s}^{(n_s)}}
$,
so by Theorem \ref{thm:borel},
$$
\tilde{\rho}_S
	\simeq
	\otimes_{s=1}^r
\bigl[
\bigl.(
  \widetilde{\sigma_{\varphi_s}^{(n_s)}} \circ \varphi_s)\bigr|_{SL_2(K)}
\bigr],
$$
where $n_s\ge 1$ and  the
$
\varphi_s: K\rarr \ov{K}
$
are field homomorphisms.   Since $SL_2(K)$ has no
non-trivial Abelian quotients, by Lemma \ref{lem:central} there is a unique
lift of
$
\tilde{\rho}_S
$
to
$
\rho|_{SL_2(K)}
$,
whence
\begin{equation}
\rho|_{SL_2(K)}
	\simeq
	\otimes_{s=1}^r
\bigl[
  \bigl.( \sigma_{\varphi_s}^{(n_s)} \circ \varphi_s)\bigr|_{SL_2(K)}
\bigr].
	\label{lift}
\end{equation}
To determine $\rho$ it then remains to study its action on the matrices
$
A_i(\alpha)
$
with $\alpha\in K^\times$.

The identity (\ref{matrix}) gives
$$
\sigma_{\varphi}^{(n)} B_{ij}(\varphi(\alpha^{\text{sgn}(j-i)}\beta))
=
\bigl[\sigma_{\varphi}^{(n)} A_i(\varphi(\alpha))\bigr]
\bigl[\sigma_{\varphi}^{(n)} B_{ij}(\varphi(\beta))\bigr]
\bigl[\sigma_{\varphi}^{(n)} A_i(\varphi(\alpha))\bigr]^{-1},
$$
whence
\begin{eqnarray*}
\lefteqn{\bigl[
  \otimes_{s=1}^r \sigma_{\varphi_s}^{(n_s)} B_{ij}(\varphi_s(\alpha^{\text{sgn}(j-i)}\beta))
\bigr]}
\\
&=&
\bigl[
  \otimes_{s=1}^r \sigma_{\varphi_s}^{(n_s)} A_i(\varphi_s(\alpha))
\bigr]
\bigl[
  \otimes_{s=1}^r \sigma_{\varphi_s}^{(n_s)} B_{ij}(\varphi_s(\beta))
\bigr]
\bigl[
  \otimes_{s=1}^r \sigma_{\varphi_s}^{(n_s)} A_i(\varphi_s(\alpha))
\bigr]^{-1}.
\end{eqnarray*}
In light of (\ref{lift}), this relation holds if we replace each
$
\otimes_{s=1}^r \sigma_{\varphi_s}^{(n_s)} B_{ij}( \varphi_s(\mbox{ -- } ))
$
by
$
\rho (B_{ij}( \mbox{ -- } ))
$.
It follows that
$$
  \bigl[\otimes_{s=1}^r \sigma_{\varphi_s}^{(n_s)} A_i(\varphi_s(\alpha))\bigr]
  \bigl[\rho(A_i(\alpha))\bigr]^{-1}
$$
commutes with every $B_{ij}(\varphi_s(\gamma))$.  As $\gamma$ runs through all
elements of
$K$, the matrices $B_{ij}(\varphi_s(\gamma))$ generate 
$
SL_2(\varphi_s(K))
$
\cite[Thm.~4.6]{artin}.  Since
$
\rho\bigr|_{SL_2(K)}
$
is absolutely irreducible, Schur's Lemma  implies that
$$
\rho(A_i(\alpha))
=
\bigl[\otimes_{s=1}^r \sigma_{\varphi_s}^{(n_s)} A_i(\varphi_s(\alpha))\bigr] f(\alpha)
$$
for some element $\alpha\in \ov{K}^\times$.  Since $\rho$ is a homomorphism,
$
f:K^\times\rarr \ov{K}^\times
$
is in fact a homomorphism of multiplicative groups.  Thus
$$
\rho = 
\bigl[
	\otimes_{s=1}^r
	\bigl( \sigma_{\varphi_s}^{(n_s)} \circ \varphi_s \bigr)
\bigr]
\otimes (f\circ\det)
$$
when restricted to $SL_2(K)$ and to the subgroup of $GL_2(K)$ generated by
$
A_i(\alpha)
$
for all the
$
\alpha\in K^\times
$.
These two subgroups together generate $GL_2(K)$, so we are done.

\end{proof}

For every $s\in SL_2(K)$, the trace of $\sigma^{(n)}(s)$ is a $\Z$-polynomial
of degree $n$ in trace$(s)$.  Write
$$
T_n( \text{trace}(s) )
$$
for this polynomial.

\begin{lem}
	\label{lem:poly}
Suppose there exists an infinite subset $\Sigma\subset SL_2(K)$ of elements
with pairwise distinct traces, such that
\begin{equation}
\prod_{i=1}^r T_{n_i}( {\rm trace}(s) )
=
\prod_{j=1}^t T_{m_j}( {\rm trace}(s) )
	\label{poly}
\end{equation}
for every $s\in\Sigma$.  Then $r=s$ and, up to reordering, $n_i=m_i$ for all
$i$.
\end{lem}

\begin{proof}
Denote by $\alpha$ and $1/\alpha$ the eigenvalues in $\ov{K}$ of $s\in\Sigma$.
Write
$
\zeta_m
$
for a fixed primitive $m$-th root of unity.  Then 
$$
T_n( \text{trace}(s) )
=
\frac{\alpha^{n+1} - (1/\alpha)^{n+1}}{\alpha - (1/\alpha)}
=
\prod_{u=1}^{n} (\alpha - \zeta_{n+1}^u (1/\alpha)).
$$
The equality (\ref{poly}) then becomes
$$
\prod_{i=1}^r \prod_{u=1}^{n_i} (\alpha - \zeta_{1+n_i}^u (1/\alpha))
=
\prod_{j=1}^t \prod_{v=1}^{m_j} (\alpha - \zeta_{1+m_j}^v (1/\alpha)).
$$
Since this holds for infinitely many distinct $\alpha\in\ov{K}$, this becomes
an equality of monic Laurent polynomials in the \textit{variable} $\alpha$ over
the field
$
L = K(\{ \zeta_{1+n_i}, \zeta_{1+m_j} \}_{i, j})
$.
Since the ring of (one variable) Laurent polynomials over $L$ has unique
factorization, we are done.

\end{proof}

\section{Odd-dimensional $l$-adic representations}

\begin{proof}[Proof of Theorem \ref{thm:main}]
Suppose $E$ is a totally real field.  We consider two cases.

\vsp

{\bf Case:}  $[E:\Q]=1$, so $\text{End}_{\ov{k}}(A)=\Z$.

Since $d$ is odd, Serre \cite[2.2.7 -- 2.2.8]{serre:onto} shows that the
image of $\rho_{A, l}$
contains $Sp_{2d}(\Z_l)$ for all but finitely many $l$.  Following the
notation in (\ref{diagram}), the restriction of the abstract homomorphism
$
\iota_\lambda
$
to
$
Sp_{2d}(\Q_l)
$
then gives rise to a non-trivial homomorphism
$
\iota'_\lambda: Sp_{2d}(\Z_l) \rarr GL_d(F_\lambda)
$
for some finite extension $F_\lambda/\Q_l$.  Since $Sp_{2d}(\Z_l)$ is compact,
the argument in
	\cite[Remark 1 on p.~1]{serre:abelian}
shows that 
$
\iota'_\lambda
$
is equivalent to a homomorphism with image in $GL_d(\oo_{F_\lambda})$.  By
Theorem \ref{thm:integral}(b), this is impossible for $l$ sufficiently large.

\vsp

{\bf Case:} $[E:\Q]=d>1$.

Suppose
$
\rho_\lambda: G_k\rarr GL_d(F_\lambda)
$
is an abstract, absolutely irreducible representation coming from $A/k$.
Following the notation in (\ref{diagram}) -- (\ref{direct}), that means
$
\rho_\lambda
$
is the composition of some
$
\tau_{l, m}: G_k\rarr GL_2(E_{l, m})
$
with an abstract homomorphism
$
\psi: GL_2(E_{l, m})\rarr G_d(\ov{\Q_l})
$.
Since $\rho_\lambda$ is absolutely irreducible, so does $\psi$.  Apply
Proposition \ref{prop:gln} and we are done.

\end{proof}

\begin{proof}[Proof of Theorem \ref{thm:indep}]
Our goal is to understand the dependency of the parameters
$r(\lambda)$,
$n_i(\lambda)$,
$f_\lambda$
and
$
\varphi_i(\lambda)
$
on $\lambda$.
Faltings shows that $\tau_{l, m}$ is semisimple, and that $E_{l, m}$ is the
commutant of
$
\tau_{l, m}(G_k)
$
in the End$(T_l(A))$.  The argument in (\cite[Thm.~4.5.4]{ribet})  then shows
that the Lie algebra of
$
\tau_{l, m}(G_k)
$
is equal to
$
\{ u\in {\mathfrak {gl}}_2(E_{l, m}):  {\rm trace}(u)\in\Q_l \}
$.
Since $\tau_{l, m}$ is equivalent to a representation whose image lies in
$
GL_2(\oo_{E_{l, m}})
$
(cf.~\cite[Remark 1 on p.~1]{serre:abelian}), it follows that
\begin{equation}
\renewcommand{\arraystretch}{1.1}
\begin{array}{ll}
\text{for every $l$, $\tau_{l, m}(G_k)$  is (conjugate to) an open}
\\
\text{subgroup in $\{ g\in GL_2(\oo_{E_{l, m}}): \det g\in \Z_l^\times \}$.}
\end{array}
	\label{image}
\renewcommand{\arraystretch}{1}
\end{equation}
In particular, there exists an infinite subset
$
\{g_i\}_i\subset G_k
$
such that, for a fixed $l$,
$
\tau_{l, m}(g_i)\in SL_2(\Z_l)\subset SL_2(\oo_{E_{l, m}})
$
and with pairwise distinct traces.  Since $\{\tau_{l, m}\}_l$ is a strictly
compatible family, that means
$
\tau_i := \text{trace}(\tau_{l, m}(g_i))
$
lies\footnote{in fact these traces lie in $\Z$,
	although we do not need it.}
in
$\Q$ and is independent of $l$.

Let $\lambda, \q$ be two distinct finite places of $F$ of residual
characteristic $l$ and $q$ ($l=q$ is allowed).  Following the notation in
Lemma \ref{lem:poly}, we have
\begin{eqnarray*}
\Phi_\lambda^{-1}
\Bigl(
  \prod_{i=1}^{r(\lambda)}  T_{n_i(\lambda)}(\tau_i)
\Bigr)
&=&
\Phi_\lambda^{-1}
\Bigl(
  \text{trace}( \rho_\lambda(g_i) )
\Bigr)
\\
&=&
\Phi_\q^{-1}
\Bigl(
  \text{trace}( \rho_\q(g_i) )
\Bigr)
	\hspace{20pt}
	\text{$\{ \rho_\lambda \}_\lambda$ is a strictly compatible family}
\\
&&
	\hspace{120pt}
	\text{and the Frobenius are dense in $G_k$}	\nonumber
\\
&=&
\Phi_\q^{-1}
\Bigl(
  \prod_{j=1}^{r(\q)}  T_{n_j(\q)}(\tau_i)
\Bigr)
\end{eqnarray*}
Since this holds for all $\tau_i$, by Lemma \ref{lem:poly} that means
$
r(\lambda) = r(\q)
$
and, up to reordering,
$
n_i(\lambda) = n_i(\q)
$
for all $i$.

It remains to study the multiplicative homomorphism
$
f_\lambda: E_l^\times\rarr\Q_l^\times
$.
Since
$
\det \tau_{l, m}
$
is the $l$-adic cyclotomic character
	\cite{ribet}, it is surjective onto $\Z_l^\times$ for
all but finitely many $l$.  Since $\Z$ is dense in $\Z_l$, we can find an
infinite subset
$
\{ h_i \}_i\subset G_k
$
and a prime $l$ such that
$
\{ \det \tau_{l, m}(h_i) \}_i = \Z - \{0\}
$.
Moreover, $d_i := \det \tau_{l, m}(h_i)$ is independent of $l$.  From
$
\text{trace}(\rho_\lambda(h_i))
=
\text{trace}(\rho_\q(h_i))
$
for all $i$ we then see that
$
f_\lambda(d_i) = f_\q(d_i)
$
for all $i$, whence
$
{f_\lambda}\bigr|_\Z = {f_\q}\bigr|_\Z
$.

\end{proof}

\section{Continunity of $l$-adic maps}

\begin{lem}
	\label{lem:lower}
For any $n>1$ there exists a constant $c_1(n)>0$, such that the degree of any
non-Abelian, complex irreducible representation of
$
SL_n(\Z/l)
$,
$
GL_n(\Z/l), Sp_{2n}(\Z/l), PSL_n(\Z/l)
$,
$PGL_n(\Z/l)$, and $PSp_{2n}(\Z/l)$ is
$
> c_1(n) l
$.
\end{lem}

\begin{proof}
Every irreducible representation of $PSL_n(\Z/l), PGL_n(\Z/l)$ and
$
PSp_{2n}(\Z/l)
$
gives rise to an irreducible representation of $SL_n(\Z/l), GL_n(\Z/l)$ and
$Sp_{2n}(\Z/l)$, respectively, so it suffices to examine the three latter
class of groups.

By Schur's Lemma, any irreducible complex representation
$
\rho: SL_n(\Z/l)\rarr GL(V)
$
takes the center of $SL_n(\Z/l)$ to the center of $GL(V)$.  Compose $\rho$
with the canonical map
$
SL(V) \rarr PGL(V)
$
and we get an irreducible \textit{projective} representation
$
PSL_n(\Z/l)\rarr PGL(V)
$.
To bound the minimal degree of non-Abelian irreducible representations of
$
SL_n(\Z/l)
$,
it then suffices to bound the minimal degree of non-Abelian irreducible
\textit{projective} representations of $PSL_n(\Z/l)$.  This is done in
	\cite{landazuri},
and the lower bound we desire follows from the main theorem there.
The same discussion applies to $Sp_{2n}$ as well.

Finally, let $\chi$ be an irreducible representation of $GL_n$.  Its
restriction to $SL_n$ remains irreducible
	\cite[Lem.~4.5]{simpson},
so $\deg \chi >1$ if $\chi$ is non-Abelian since $SL_n$ is its own
commutator subgroup.  The Lemma for $GL_n$ is then reduced to the case of
$SL_n$.  This completes the proof of the Lemma.

\end{proof}

\begin{lem}
	\label{lem:nomap}
Fix integers $m, N\ge 1$ and $n\ge 2$. Then there exists a constant
$
c_2(m, n, N)>0
$
such that, for any prime $l>c_2(m, n, N)$, no subgroup of $GL_N(\C)$ is
isomorphic to
$
SL_n(\Z/l^m), GL_n(\Z/l^m), Sp_{2n}(\Z/l^m), PSL_n(\Z/l^m), PGL_n(\Z/l^m)
$,
or $PSp_{2n}(\Z/l^m)$.
\end{lem}

\begin{proof}
Suppose there exists an injective map
$
\phi: SL_n(\Z/l^m)\rarr GL_N(\C)
$.
We can assume that $\phi$ is irreducible.
Denote by $S$ the kernel of the reduction map
$
SL_n(\Z/l^m)\rarr SL_n(\Z/l)
$.
Every irreducible representation of a nilpotent group $P$ is induced from a
degree one character of some subgroup of $P$
	\cite[Thm.~52.1]{curtis},
so the degree of every irreducible representation of the $l$-group $S$ is
either $1$ or  is divisible by
$l$.  But 
$
\deg (\phi\bigr|_S) = N
$,
so for $l>N$ every irreducible representation of $S$ appearing in
$
\phi\bigr|_S
$
has degree $1$.  Thus $\phi(S)$ is Abelian.  But $S$ is non-Abelian if
$m>2$, so $m\le 2$.  Lemma \ref{lem:lower} covers the case $m=1$, so
it remains to consider the case $m=2$.  We first introduce some notation.

\vsp

Let $H$ be a
normal subgroup of a finite group $G$.  For any character $\psi$ and any
element $g\in G$, denote by
$
\psi^g
$
the \textit{conjugate} of $\psi$ by $g$:
\begin{equation}
\psi^g(h) := \psi(g^{-1}hg)	\hspace{10pt}\text{for all $h\in H$.}
	\label{conj}
\end{equation}

\vsp

We now resume the proof of the Lemma.  Take $m=2$, so $S$ is Abelian.
Clifford's theory \cite[$\S 49$]{curtis} gives a decomposition
\begin{equation}
\phi\bigr|_S \simeq (\psi^{s_1} \oplus \cdots \oplus \psi^{s_r})^e,
	\label{cliff}
\end{equation}
where $e, r\ge 1$ are integers,
$
s_i\in SL_n(\Z/l^2)
$,
and $\psi$ is an irreducible representation of $S$.  Moreover,
$
SL_n(\Z/l^2)
$
permutes
the $\psi^{g_i}$ transitively via (\ref{conj}).  We claim that this action on
the $\psi^{g_i}$, when restricted to $S$, is trivial.
To see this, note that for any $h, s\in S$,
\begin{eqnarray*}
(\psi^{g_i})^{s}(h)
&=&
\psi(s^{-1} g_i^{-1} h g_i s)
\\
&=&
\psi(g_i^{-1} \sigma^{-1} h \sigma g_i)
	\hspace{10pt}
	\text{for some $\sigma\in S$, since $S$ is normal}
\\
&=&
\psi(g_i^{-1} h g_i)
	\hspace{36pt}
	\text{since $S$ is Abelian}
\\
&=&
\psi^{g_i}(h).
\end{eqnarray*}
This verifies the claim.  As a result,
$
SL_n(\Z/l^2)/S \simeq SL_n(\Z/l)
$
permutes the $\psi^{g_i}$ transitively.  The left side of (\ref{cliff}) has
degree $N$, so $r\le N$.
Thus $SL_n(\Z/l)$ permutes transitively a set of size $r\le N$.  We have
two cases:

\begin{itemize}
\item
$r=1$:
\\
Then
$
\phi\bigr|_S
$
is the direct sum of copies of a fixed irreducible representation $\psi$ of
$S$.  Since $S$ is Abelian, that means $\deg \psi = 1$, whence
$
\phi\bigr|_S
$
has cyclic image.  This is impossible since $S$ is non-cyclic and $\phi$ is
injective.
\item
$1<r\le N$:
\\
Then $SL_n(\Z/l)$ contains a proper subgroup of index $\le N$.  That means
$
SL_n(\Z/l)
$
has a non-trivial permutation representation, and hence a non-trivial
irreducible
representation, of degree $\le N$.  This is impossible if $l$ is large
enough, by Lemma \ref{lem:lower}
(note that $SL_n(\Z/l)$ has no non-trivial degree one character if
$
(n, l) \not= (2,2)
$).
\end{itemize}
This completes the proof of the Lemma for $SL_n$.  The same argument applies
to $PSL_n, Sp_{2n}$ and $PSp_{2n}$.  Finally, if $GL_N(\C)$ contains a subgroup
isomorphic to $GL_n(\Z/l^m)$ or $PGL_n(\Z/l^m)$, then it contains one
isomorphic to $SL_n(\Z/l^m)$ or $PSL_n(\Z/l^m)$, so we are done.

\end{proof}

\if 3\
{
It remains to show that 
$
\ker\phi'
$
contains $S_m\simeq (\Z/l)^{\oplus 3}$.  Suppose otherwise.  Then
$
\#(\ker\phi'\cap S_m) = 1, l
$
or $l^2$.   If
$
\ker\phi'\cap S_m
$
is trivial,  then \cite[Prop.~2.1]{dennin} implies
that
$
\ker\phi'\cap S = \ker\phi'
$
is trivial, which is impossible.   If
$
\#(\ker\phi'\cap S_m) = l
$,
then \cite[Prop.~2.2]{dennin} implies that $\ker\phi'$ is cyclic, whence
$
\text{im}\,\phi' \simeq S/\ker\phi'
$
is not cyclic, contradicting the fact that $\phi'$ is an Abelian character.
Finally, suppose
$
\#(\ker\phi'\cap S_m) = l^2
$.
By \cite[Lem.~4]{mcquillan}
$
\ker\phi'\cap S_m
$
is not a normal subgroup of $SL_2(\Z/l^m)$, and hence not a normal subgroup of
$S$, contradicting the fact that
$
\ker\phi'
$
is normal in $S$.  This shows that
$
\phi'\bigr|_{S_m}
$
is trivial, as desired.
}
\fi

\begin{prop}
	\label{prop:local}
Let $A$ be a local ring with maximal ideal $\lambda$.  Suppose
$
\text{char}(A/\lambda)>3
$.
Then a subgroup $H\subset GL_N(A)$ is normal if and only if there exists an
ideal $\aaa\subset A$ such that the image of $H$ in $GL_N(A/\aaa)$ consists
of the scalar matrices.
The same holds for
$
SL_N, Sp_{2N}, PGL_N, PSL_N
$
and $PSp_{2N}$.
\end{prop}

\begin{proof}
For $GL_N(A)$ and $SL_N(A)$ see \cite[p.~84 and p.~245]{mcdonald}.  For
$Sp_{2d}(A)$  see \cite[p.~210]{mcdonald}.

\end{proof}

\begin{proof}[Proof of Theorem \ref{thm:integral}]
Again we give the argument for $SL_n$.  We can assume that $l>3$.
Denote by $\lambda$ and $\lambda'$ the maximal ideal of $\oo$ and $\oo'$,
respectively.  

\vsp

(a)
We can assume that $\varphi$ is not trivial.
For $i> 0$, denote by $G_i$ the kernel of the reduction map
$
GL_N(\oo) \rarr GL_N(\oo/\lambda^i)
$
Every quotient $G_i/G_{i+1}$ is a finite Abelian $l$-group, whence the same
holds for
$
\varphi^{-1}(G_i)/\varphi^{-1}(G_{i+1})
$.

\vsp

\begin{lem}
	\label{lem:claim}
Every $\varphi^{-1}(G_i)$ has finite index in $SL_n(\oo')$.
\end{lem}

\begin{proof}
Every $G_i$ is normal in $GL_N(\oo)$, so every
$
\varphi^{-1}(G_i)
$
is normal in $SL_n(\oo')$.  In light of Proposition \ref{prop:local},
it suffices to show that none of the 
$
\varphi^{-1}(G_i)
$
is contained in the center of $SL_n(\oo')$.  Suppose otherwise.  Then
$
SL_n(\oo')/\varphi^{-1}(G_i)
$
is infinite and injects into 
$
GL_N(\oo)/G_i \simeq GL_N(\oo/\lambda^i)
$,
which is finite.  This is a contradiction.

\end{proof}

The $G_i$ form a basis of open
neighborhoods of the identity element in $GL_N(\oo)$.  Lemma \ref{lem:claim}
then implies that
$
\varphi^{-1}(U)
$
has finite index in $SL_n(\oo')$ for every sufficiently small open set
$
U\subset GL_N(\oo)
$
containing the identity; and hence for all open sets $U\subset GL_N(\oo)$.
That means $\varphi$ is continuous.

Finally, suppose $\ker\varphi$ is not contained in the subgroup of scalar
matrices.   In light of
Proposition \ref{prop:local}, for some $m\ge 1$ we then get an injective
map of
$
SL_n(\oo'/{\lambda'\,}^m)
$
or
$
PSL_n(\oo'/{\lambda'\,}^m)
$
into $GL_N(\oo)$.  Restrict this map to the subgroup
$
SL_n(\Z/l^m)
$
or
$
PSL_n(\Z/l^m)
$
and note that we have a (discontinuous) inclusion
$
GL_N(\oo)\subset GL_N(\C)
$,
we get a contradiction, by Lemma \ref{lem:nomap}.  This completes the proof
of Part (a).

\vsp

(b) We claim that if such $\varphi$ exists, then we can find an injective map
of abstract
groups
\begin{eqnarray*}
SL_n(\F_l)  \hookrightarrow GL_N(\ov{\F}_l)
	&
	\text{or}
	&
PSL_n(\F_l) \hookrightarrow GL_N(\ov{\F}_l).
\end{eqnarray*}
By the work of Steinberg \cite{steinberg} it then follows that $N\ge n$.

\vsp

For every integer $i\ge 1$, denote by $G_i$ the kernel of the reduction map
$
GL_N(\oo)\rarr GL_N(\oo/\lambda^i)
$.
Every quotient $G_i/G_{i+1}$ is a finite Abelian $l$-group, whence the same
holds for the
$
{\varphi}^{-1}(G_i)/ {\varphi}^{-1}(G_{i+1})
$.
We claim that
\begin{equation}
{\varphi}^{-1}(G_1) \not= SL_n(\oo').
	\label{otherwise}
\end{equation}
Suppose otherwise.  Since $\text{ker}(\varphi)$ is a proper subgroup of
$
SL_n(\oo')
$
and since
$
\cap_i G_i = \{ I \}
$,
if (\ref{otherwise}) is false then there exists a smallest integer $j\ge 1$
such that
$
{\varphi}^{-1}(G_{j+1})\subsetneq  {\varphi}^{-1}(G_j)
$;
necessarily
$
{\varphi}^{-1}(G_j) = SL_n(\oo')
$.
On the other hand, every $\varphi^{-1}(G_i)$ is normal in
$
SL_n(\oo')
$.
Since $l>3$, Proposition \ref{prop:local} implies that
$$
\varphi^{-1}(G_{j+1})
=
\ker( 
\raisebox{2.5pt}{$\begin{CD}
SL_n(\oo') @>{\text{mod }{\lambda'\,}^m}>> SL_n(\oo'/{\lambda'\,}^m))
\end{CD}$}
)
$$
or
$$
\varphi^{-1}(G_{j+1})
=
\ker( 
\raisebox{2.5pt}{$\begin{CD}
SL_n(\oo') @>{\text{mod }{\lambda'\,}^m}>> PSL_n(\oo'/{\lambda'\,}^m))
\end{CD}$}
)
$$
for some $m\ge 1$.  Then
$
{\varphi}^{-1}(G_j) / {\varphi}^{-1}(G_{j+1})
=
SL_n(\oo') / {\varphi}^{-1}(G_{j+1})
\simeq
SL_n(\oo'/{\lambda'\,}^m)
$
or 
$
PSL_n(\oo'/{\lambda'\,}^m)
$,
\if 3\
{
\cite[p.~212]{mcdonald} the proper normal subgroups of
$
SL_n(\Z_l)
$
are given by
\begin{equation}
  \renewcommand{\arraystretch}{1.2}
  \begin{array}{lllllll}
	\{ x\in SL_n(\Z_l): x\equiv\mymod{I}{\mathfrak a} \}
	&
		\text{and}
	&
	\{ x\in SL_n(\Z_l): x\equiv\mymod{\pm I}{\mathfrak a} \}
  \end{array}
  \renewcommand{\arraystretch}{1}
	\label{normal}
\end{equation}
as we run through all ideals $\mathfrak a\subsetneq \Z_l$  (we recover the
subgroups
$
\{ I \}
$
and
$
\{ \pm I \}
$
by setting $\mathfrak a = (0)$).  Then
}
\fi
which is not a finite Abelian $l$-group.  This is a contradiction, so
(\ref{otherwise}) must hold, and we get an injective map
$$
SL_n(\oo') / {\varphi}^{-1}(G_1)
	\hookrightarrow
GL_N(\oo_L)/G_1
	\simeq
GL_N(\oo_L/\lambda).
$$
By Proposition \ref{prop:local}, we get an
injective
map
\begin{eqnarray}
SL_n(\oo'/{\lambda'\,}^m) \hookrightarrow GL_N(\oo_L/\lambda)
	&\text{or}&
PSL_n(\oo'/{\lambda'\,}^m) \hookrightarrow GL_N(\oo_L/\lambda).
	\label{hook2}
\end{eqnarray}
for some $m\ge 1$.  Every element of
$
GL_N(\oo_L/\lambda)
$
either is semisimple or contains a non-trivial Jordan block
	(over $\ov{\F}_l$).
In the first case the
element has order prime to $l$.  In the second case, if $l>N$ then every
non-semisimple element has order divisible by $l$ but not by $l^2$.
On the other hand, since $l>2$, both
$
SL_n(\oo'/{\lambda'\,}^m)
$
and
$
PSL_n(\oo'/{\lambda'\,}^m)
$
contain elements of order $l^m$ (coming from transvections).
Since the maps in (\ref{hook2}) are injective, that means $m=1$.
This completes the proof of Part (b).

\end{proof}

\if 3\
{
\vsp

We now tackle (ii).
Since $l\not=3$,  we can choose an odd prime $q\not=l$ that
divides the order of $SL_2(\F_l) = Sp_2(\F_l)$. Then $SL_n(\F_l)$ has a
subgroup $A$ which is isomorphic to $(\Z/q)^d$.  So does $PSL_n(\F_l)$.
But in $GL_N(\ov{\F}_l)$ there
is, up to conjugation, exactly one subgroup isomorphic to $(\Z/q)^d$, and it
intersects the center of $GL_N(\ov{\F}_l)$ in a group of order $q$. But
neither $SL_n(\F_l)$ nor $PSL_n(\F_l)$ has a central subgroup of order
$q>2$, a
contradiction.  This gives (ii), and hence completes the proof of the Lemma.
}
\fi

\section*{Acknowledgement}
My interest in odd-dimensional Galois representations was inspired by a
lecture of Ash
on the analog of Serre's conjecture for $GL(3)$ modular representations.
I am indebted to Humphreys for his help with the work of
Borel-Tits,  and to Ribet for his help concerning the image of
$\lambda$-adic representations.  I would like to thank 
Dettweiler for pointing out a mistake in an earlier draft.
I also benefited from a useful conversation with Rosen.

\bibliographystyle{amsalpha}

\end{document}